\documentclass[a4paper,11pt]{amsart}
\usepackage{graphicx}
 \usepackage{mathptmx}
\usepackage{amsmath}
\usepackage{amssymb}
\usepackage{enumitem}
\usepackage{xcolor}

\newmuskip\pFqmuskip

\newcommand*\pFq[6][8]{%
  \begingroup 
  \pFqmuskip=#1mu\relax
  \mathcode`=\string"8000
  \begingroup\lccode`\~=`\,
  \lowercase{\endgroup\let~}\pFqcomma
  F^{#2}_{#3}{\left(\genfrac..{0pt}{}{#4}{#5}\bigg|#6\right)}%
  \endgroup
}
\newcommand{\pFqcomma}{\mskip\pFqmuskip}

\newtheorem{theorem}{Theorem}[section]

\newtheorem{corollary}[theorem]{Corollary}

\begin{document}

\title{Degenerate hypergeometric functions and degenerate hypergeometric numbers of order $\bf{p}$}

\author{Taekyun  Kim}
\address{Department of Mathematics, Kwangwoon University, Seoul 139-701, Republic of Korea}
\email{tkkim@kw.ac.kr}

\author{Dae San  Kim }
\address{Department of Mathematics, Sogang University, Seoul 121-742, Republic of Korea}
\email{dskim@sogang.ac.kr}
\author{Hyunseok Lee}
\address{Department of Mathematics, Kwangwoon University, Seoul 139-701, Republic of Korea}
\email{luciasconstant@gmail.com}


\maketitle

\begin{abstract}
Recently, Simsek studied certain finite sums involving powers of binomials coefficients which are called generalized $p$-th order Franel numbers and can be represented in terms of hypergeometric functions. Then, among other things, he showed that particular cases of these numbers are connected with many known special numbers and polynomials. The purpose of this paper is to investigate degenerate versions of those numbers.
In more detail, we introduce degenerate generalized hypergeometric functions and study degenerate hypergeometric numbers of order $p$. These numbers involve powers of $\lambda$-binomial coefficients and $\lambda$-falling sequence, and can be represented by means of the degenerate generalized hypergeometric functions. We will derive some explicit expressions and combinatorial identities for those numbers. We also consider several related special numbers like $\lambda$-hypergeometric numbers of order $p$ and Apostol type $\lambda$-hypergeometric numbers of order $p$, of which the latter reduce in a limiting case to the generalized $p$-th order Franel numbers.
\end{abstract}

\pagestyle{myheadings}
\markboth{\centerline{\scriptsize T. Kim, D.S. Kim, H. Lee}}
          {\centerline{\scriptsize Degenerate hypergeometric functions and degenerate hypergeometric numbers of order $\bf{p}$}}

\section{Introduction}
Recently, Simsek \cite{17} studied certain finite sums involving powers of binomials coefficients which are called generalized $p$-th order Franel numbers and can be represented in terms of hypergeometric functions. Then, among other things, he showed that particular cases of these numbers are connected with many known special numbers and polynomials which include  Bernoulli numbers, Euler numbers, Changhee numbers, Daehee numbers, Stirling numbers of the first kind, Catalan numbers and Legendre polynomials.\\
\indent In recent years, many mathematicians have drawn their attention in studying various degenerate versions of some special numbers and polynomials \cite{3,6,11,12,14}. The idea of investigating degenerate versions of some special numbers and polynomials originated from Carlitz's papers \cite{1,2}. Indeed, he introduced the degenerate Bernoulli and Euler polynomials and numbers, and investigated some arithmetic and combinatorial aspects of them.
Here we mention in passing that the degenerate Bernoulli polynomials were later rediscovered by Ustinov under the name of Korobov polynomoals of the second.
Two of the present authors, their colleagues and some other people have studied quite a few degenerate versions of special numbers and polynomials with their interest in not only combinatorial and arithmetic properties and but also in differential equations and certain symmetric identities [8,14 and references therein]. It is worth noting that this idea of considering degenerate versions of some special polynomials and numbres is not only limited to polynomials but also can be extended to transcendental functions like gamma functions \cite{9,10}.  We believe that studying some degenerate versions of special polynomials and numbers is very fruitful and promising area of research in which many things remain yet to be uncovered. \\
\indent Motivated by Simsek's paper [17], we would like to investigate degenerate versions of the generalized $p$-th order Franel numbers. In more detail, we introduce degenerate generalized hypergeometric functions and study degenerate hypergeometric numbers of order $p$. These numbers involve powers of $\lambda$-binomial coefficients and $\lambda$-falling sequence, and can be represented by means of the degenerate generalized hypergeometric functions. We also consider several related special numbers like $\lambda$-hypergeometric numbers of order $p$ and Apostol type $\lambda$-hypergeometric numbers of order $p$ of which the latter reduces in a limiting case to the generalized $p$-th order Franel numbers. \\
\indent For the rest of this section, we will fix some notations and recall some known results that are needed throughout this paper.

\vspace{0.1in}

For $\lambda\in\mathbb{R}$, the degenerate exponential function is defined as
\begin{equation}
	\label{1} e_{\lambda}^{x}(t)=(1+\lambda t)^{\frac{x}{\lambda}},\quad e_{\lambda}(t)=(1+\lambda t)^{\frac{1}{\lambda}}=e_{\lambda}^{1}(t),\ (\mathrm{see},\ [1,2,9,10]).
\end{equation}
From \eqref{1}, we note that
\begin{equation}
	\label{2} e_{\lambda}^{x}(t)=\sum_{n=0}^{\infty}(x)_{n,\lambda}\frac{t^{n}}{n!},\quad (\mathrm{see}\ [3,8,9,10,11]),
\end{equation}
where $(x)_{n,\lambda}$ are the $\lambda$-falling sequence given by
\begin{equation}
	\label{3} (x)_{0,\lambda}=1,\quad (x)_{n,\lambda}=x(x-\lambda)\cdots(x-(n-1)\lambda),\quad (n\ge 1).
\end{equation}
\indent In [11], the degenerate Stirling numbers of the second kind are defined by
\begin{equation}
	\label{4} \frac{1}{k!}(e_{\lambda}(t)-1)^{k}=\sum_{n=k}^{\infty}S_{2,\lambda}(n,k)\frac{t^{n}}{n!},\quad (k\ge 0).
\end{equation}
Let
\begin{equation}
	\label{5} (x)_{0}=1,\quad (x)_{n}=x(x-1)(x-2)\cdots(x-(n-1)),\quad (n\ge 1).
\end{equation}
Then $\displaystyle\lim_{\lambda\rightarrow 0}S_{2,\lambda}(n,k)=S_{2}(n,k)\displaystyle$, where $S_{2}(n,k)$ are the ordinary Stirling numbers of the second kind given by
\begin{equation}
	\label{6} x^{n}=\sum_{l=0}^{n}S_{2}(n,l)(x)_{l},\quad (n\ge 0),\quad (\mathrm{see}\ [1-21]).
\end{equation}

The Stirling numbers of the first kind are defined as
\begin{equation}
	\label{7} \frac{1}{k!}\big(\log(1+t)\big)^{k}=\sum_{n=k}^{\infty}S_{1}(n,k)\frac{t^{n}}{n!},\quad (k\ge 0),\ (\mathrm{see}\ [12,13]).
\end{equation}
Thus, by \eqref{7}, we get
\begin{equation}
	\label{8} (x)_{n}=\sum_{l=0}^{n}S_{1}(n,l)x^{l},\quad (n\ge 0).
\end{equation}
In view of \eqref{4}, the degenerate Stirling numbers of the first kind are defined by
\begin{equation}
	\label{9} \frac{1}{k!}\bigg(\frac{t^{\lambda}-1}{\lambda}\bigg)^{k}=\sum_{n=k}^{\infty}S_{1,\lambda}(n,k)\frac{t^{n}}{n!},\quad (\mathrm{see}\ [3,11]).
	\end{equation}
Note that $\displaystyle\lim_{\lambda\rightarrow 0}S_{1,\lambda}(n,k)=S_{1}(n,k),\ (n,k\ge 0)\displaystyle$.\\

\indent As is well known, the generalized hypergeometric function $F^{(p,q)}$ is defined by  	
\begin{align}\label{10}
\pFq{(p,q)}{}{a_{1},a_{2},\dots,a_{p}}{b_{1},b_{2},\dots,b_{q}}{x}
=\sum_{k=0}^{\infty}\frac{\langle a_{1}\rangle_{k}\cdots\langle a_{p}\rangle_{k}}{\langle b_{1}\rangle_{k} \langle b_{2}\rangle_{k}\cdots \langle b_{q}\rangle_{k}}\frac{x^{k}}{k!},
\end{align}
where $\langle a\rangle_{k}=a(a+1)\cdots(a+(k-1)),\, (k\ge 1)$,\, $ \langle a\rangle_{0}=1$, (see [15, 20]). \\
For example,
\begin{displaymath}
	\pFq{(2,1)}{}{1,b}{b}{-1}=\frac{1}{2},\ \ \pFq{(2,1)}{}{2,b}{b}{-1}=\frac{1}{4},\ \ \pFq{(2,1)}{}{3,b}{b}{-1}=\frac{1}{8},\cdots.
\end{displaymath}
The Gauss summation theorem is given by
\begin{equation}
	\label{11} \pFq{(2,1)}{}{a,b}{c}{1}=\frac{\Gamma(c)\Gamma(c-a-b)}{\Gamma(c-a)\Gamma(c-b)},
\end{equation}
where $R(c)>R(b)>0$, $\ R(c-a-b)>0$, $\ R(c)>R(a)>0$.\\
\
 From \eqref{11}, we note that
\begin{equation}
	\label{12} \pFq{(2,1)}{}{a,b}{c}{z}=\frac{\Gamma(c)}{\Gamma(b)\Gamma(c-b)}\int_{0}^{1}t^{b-1}(1-t)^{c-b-1}(1-tz)^{-a}dt
\end{equation}
where $R(c)>R(b)>0$. \\

\vspace{0.1in}

The following are well known identities related to the binomial coefficients:
\begin{align}	
\label{13} \sum_{k=0}^{n}\binom{n}{k}&=2^{n},\quad (n\ge 0),\\
\label{14} \sum_{k=0}^{n}(-1)^{k}\binom{n}{k}&=0,\quad (n\ne 0,\ n\in\mathbb{N}), \\
\label{15} \sum_{k=0}^{n}\binom{n}{k}^{2}&=\frac{(2n)!}{(n!)^{2}}=\binom{2n}{n}, \\
\label{16}\sum_{k=0}^{n}(-1)^{k}\binom{n}{k}^{3}&=\left\{\begin{array}{cc} 0, & \textrm{if $n$ is odd}\\ \frac{(-1)^{n/2}n!}{\big((\frac{n}{2})!\big)^{2}} & \textrm{otherwise.}\end{array}\right. \\
\label{17}\sum_{k=-n}^{n}(-1)^{k}\binom{n+b}{n+k}\binom{n+c}{c+k}\binom{b+c}{b+k}&=\frac{\Gamma(b+c+n+1)}{n!\Gamma(b+1)\Gamma(c+1)}.
\end{align}

\vspace{0.1in}

\section{Sums of powers of $\lambda$-binomial coefficients}
The $\lambda$-binomial coefficients are defined as
\begin{equation}
	\label{18}\binom{x}{n}_{\lambda}=\frac{(x)_{n,\lambda}}{n!}=\frac{x(x-\lambda)\cdots(x-(n-1)\lambda)}{n!},\quad (n\ge 1),\quad \binom{x}{0}_{\lambda}=1,\quad (\lambda\in\mathbb{R}).
\end{equation}
From \eqref{18}, we easily get
\begin{equation}
	\label{19}\binom{x+y}{n}_{\lambda}=\sum_{l=0}^{n}\binom{x}{l}_{\lambda}\binom{y}{n-l}_{\lambda},\quad (n\ge 0).
\end{equation}
\indent By \eqref{1}, we easily get
\begin{equation}
	\label{20}B_{\lambda}^{*}(n,k)=\frac{d^{n}}{dt^{n}}\big(e_{\lambda}(t)+1\big)^{k}\bigg|_{t=0}=\sum_{j=1}^{k}\binom{k}{j}(j)_{n,\lambda},
\end{equation}
where $n$ and $k$ are positive integers. \\
Note that $\displaystyle\lim_{\lambda\rightarrow 0}B_{\lambda}^{*}(n,k)=B(n,k)$, where $B(n,k)$ are defined by Golombek and given by
\begin{displaymath}
	B(n,k)=\sum_{j=1}^{k}\binom{k}{j}j^{n}\quad (\mathrm{see}\ [4,5]).
\end{displaymath}
\indent Now, we define the {\it{degenerate hypergeometric function}} as
\begin{equation}
\label{21}	\pFq{}{\lambda}{a,b}{c}{z}=\sum_{n=0}^{\infty}\frac{\langle a\rangle_{n,\lambda}\langle b\rangle_{n,\lambda}}{\langle c\rangle_{n,\lambda}}\frac{z^{n}}{n!},
\end{equation}
where $\langle a\rangle_{n,\lambda}=a(a+\lambda)\cdots(a+(n-1)\lambda),\ (n\ge 1)$, $\langle a\rangle_{0,\lambda}=1$.\\
From \eqref{21}, we note that
\begin{align}
	\label{22} \pFq{}{\lambda}{-n,-n}{\lambda}{\lambda e_{\lambda}(t)}&=\sum_{k=0}^{\infty}\frac{\langle -n\rangle_{k,\lambda}\langle -n\rangle_{k,\lambda}\lambda^{k}}{\langle\lambda\rangle_{k,\lambda}}\frac{e_{\lambda}^{k}(t)}{k!}=\sum_{k=0}^{\infty}\frac{(n)_{k,\lambda}(n)_{k,\lambda}}{k!}\frac{e_{\lambda}^{k}(t)}{k!}\\
\nonumber &=\sum_{k=0}^{\infty}\binom{n}{k}_{\lambda}^{2}e^{k}_{\lambda}(t)=\sum_{m=0}^{\infty}\sum_{k=0}^{\infty}\binom{n}{k}_{\lambda}^{2}(k)_{m,\lambda}\frac{t^{m}}{m!},
\end{align}
where $n$ is a nonnegative integer. \\
\indent Let us define
\begin{equation}
	\label{23} \pFq{}{\lambda}{-n,-n}{\lambda}{\lambda e_{\lambda}(t)}=\sum_{m=0}^{\infty}Q_{\lambda}(m,2)\frac{t^{m}}{m!}.
\end{equation}
Therefore, by (22) and (23), we obtain the following theorem.
\begin{theorem}
	For $m\ge 0$, we have
	\begin{displaymath}
		Q_{\lambda}(m,2)=\sum_{k=0}^{\infty}\binom{n}{k}_{\lambda}^{2}(k)_{m,\lambda}.
	\end{displaymath}
\end{theorem}
\noindent We note that $\displaystyle\lim_{\lambda\rightarrow 0}Q_{\lambda}(m,2)=\sum_{k=0}^{n}\binom{n}{k}^{2}k^{m}=Q(m,2)\displaystyle$, which was introduced by Golombek and Marburg (see [4,5]). \\
\indent We observe that
\begin{displaymath}
	\pFq{}{\lambda}{-a}{\times}{-z}=\sum_{k=0}^{\infty}\langle -a\rangle_{k,\lambda}\frac{(-z)^{k}}{k!}=\sum_{k=0}^{\infty}\frac{(a)_{k,\lambda}}{k!}z^{k}=\sum_{k=0}^{\infty}\binom{a}{k}_{\lambda}z^{k}=e_{\lambda}^{a}(z).
\end{displaymath}
For $n\in\mathbb{N}$, let
\begin{equation}
\label{24} \pFq{}{\lambda}{-n}{\times}{-e_{\lambda}(t)}=(1+\lambda e_{\lambda}(t))^{\frac{n}{\lambda}}=\sum_{m=0}^{\infty}H_{\lambda}(n,m)\frac{t^{m}}{m!}.
\end{equation}
On the one hand, we have
\begin{align}	
\label{25}(1+\lambda e_{\lambda}(t))^{\frac{n}{\lambda}}=\sum_{k=0}^{\infty}\binom{n}{k}_{\lambda}e^{k}_{\lambda}(t)=\sum_{m=0}^{\infty}\bigg(\sum_{k=0}^{\infty}\binom{n}{k}_{\lambda}(k)_{m,\lambda}\bigg)\frac{t^{m}}{m!}.\end{align}
On the other hand, we get
\begin{align}
	\label{26} (1+\lambda e_{\lambda}(t))^{\frac{n}{\lambda}}& =\big(1+\lambda+\lambda(e_{\lambda}(t)-1)\big)^{\frac{n}{\lambda}}\\
	\nonumber &=(1+\lambda)^{\frac{n}{\lambda}}\bigg(1+\frac{\lambda}{1+\lambda}(e_{\lambda}(t)-1)\bigg)^{\frac{n}{\lambda}}\\
	\nonumber &= (1+\lambda)^{\frac{n}{\lambda}}\sum_{k=0}^{\infty}\binom{n}{k}_{\lambda}\bigg(\frac{1}{1+k}\bigg)^{k}(e_{\lambda}(t)-1)^{k}\\
	\nonumber &=(1+\lambda)^{\frac{n}{\lambda}}\sum_{k=0}^{\infty}(n)_{k,\lambda}\bigg(\frac{1}{1+\lambda}\bigg)^{k}\frac{1}{k!}(e_{\lambda}(t)-1)^{k}\\
	\nonumber &=\sum_{m=0}^{\infty}\bigg((1+\lambda)^{\frac{n}{\lambda}}\sum_{k=0}^{m}(n)_{k,\lambda}\bigg(\frac{1}{1+\lambda}\bigg)^{k}S_{2,\lambda}(m,k)\bigg)\frac{t^{m}}{m!}.
\end{align}
From \eqref{24}, \eqref{25} and \eqref{26}, we obtain the following theorem.
\begin{theorem}
	For $n\in\mathbb{N}$ and $m\in\mathbb{N}\cup\{0\}$, we have
	\begin{displaymath}
		H_{\lambda}(n,m)=(1+\lambda)^{\frac{n}{\lambda}}\sum_{k=0}^{m}(n)_{k,\lambda}\bigg(\frac{1}{1+\lambda}\bigg)^{k}S_{2,\lambda}(m,k)=\sum_{k=0}^{\infty}\binom{n}{k}_{\lambda}(k)_{m,\lambda}.
	\end{displaymath}
\end{theorem}
As is well known, the degenerate Bell polynomials are defined by
\begin{equation}
	\label{28} e_{\lambda}^{x}(e_{\lambda}(t)-1)=\sum_{n=0}^{\infty}\mathrm{Bel}_{n,\lambda}(x)\frac{t^{n}}{n!},\quad (\mathrm{see}\ [14]).
\end{equation}
By \eqref{28}, we easily get
\begin{equation}
	\label{29} \mathrm{Bel}_{n,\lambda}(x)=\sum_{k=0}^{n}(x)_{k,\lambda}S_{2,\lambda}(n,k),\quad (n\ge 0).
\end{equation}
Now, we define the {\it{degenerate bivariate Bell polynomials}} by
\begin{equation}
	\label{30} e_{\lambda}^{x}\big(y(e_{\lambda}(t)-1)\big)=\sum_{n=0}^{\infty}\mathrm{Bel}_{n,\lambda}(x,y)\frac{t^{n}}{n!}.
\end{equation}
Thus, by (30), we get
\begin{equation}
	\label{31} \mathrm{Bel}_{n,\lambda}(x,y)=\sum_{k=0}^{n}(x)_{k,\lambda}y^{k}S_{2,\lambda}(n,k),\quad (n\ge 0).
\end{equation}
From Theorem 2.2 and \eqref{31}, we obtain the following corollary.
\begin{corollary}
For $n\in\mathbb{N}$ and $m\in\mathbb{N}\cup\{0\}$, we have
\begin{displaymath}
H_{\lambda}(n,m)=(1+\lambda)^{\frac{n}{\lambda}}\mathrm{Bel}_{m,\lambda}\bigg(n,\frac{1}{1+\lambda}\bigg).	
\end{displaymath}
\end{corollary}
Note that
\begin{displaymath}
\lim_{\lambda\rightarrow 0}H_{\lambda}(n,m)=\sum_{k=0}^{\infty}\frac{n^k}{k!}k^{m}=e^{n}\sum_{k=0}^{m}n^{k}S_{2}(m,k).
\end{displaymath}
We observe that
\begin{align}
\nonumber H_{\lambda}(n,1)&=\sum_{k=0}^{\infty}\binom{n}{k}_{\lambda}(k)_{1,\lambda}=\sum_{k=1}^{\infty}\frac{(n)_{k,\lambda}}{(k-1)!}\\\label{32}	&=\sum_{k=0}^{\infty}\binom{n}{k}_{\lambda}(n-k\lambda)=n\sum_{k=0}^{\infty}\binom{n}{k}_{\lambda}-\lambda\sum_{k=0}^{\infty}\binom{n}{k}_{\lambda}k \\
\nonumber &=n\sum_{k=0}^{\infty}\binom{n}{k}_{\lambda}-\lambda H_{\lambda}(n,1)
\end{align}
Thus, by \eqref{32}, we get
\begin{equation}
	\label{33} (1+\lambda)H_{\lambda}(n,1)=n\sum_{k=0}^{\infty}\binom{n}{k}_{\lambda}=n(1+\lambda)^{\frac{n}{\lambda}}
\end{equation}
From \eqref{33}, we have
\begin{displaymath}
H_{\lambda}(n,1)=n(1+\lambda)^{\frac{n}{\lambda}-1},\,\,\big(\lambda \neq -1\big).
\end{displaymath}
\indent For $m=2$, we have
\begin{align}
\label{34} \sum_{k=0}^{\infty}\binom{n}{k}_{\lambda}(k)_{2,\lambda}&=\sum_{k=0}^{\infty}\frac{(n)_{k,\lambda}}{k!}k(k-\lambda)=\sum_{k=1}^{\infty}\frac{(n)_{k,\lambda}}{(k-1)!}(k-\lambda)\\ \nonumber
&=\sum_{k=1}^{\infty}\frac{(n)_{k,\lambda}}{(k-1)!}(k-1+1-\lambda)
=\sum_{k=2}^{\infty}\frac{(n)_{k,\lambda}}{(k-2)!}+(1-\lambda)\sum_{k=1}^{\infty}\frac{(n)_{k,\lambda}}{(k-1)!}\\
\nonumber &=\sum_{k=2}^{\infty}\frac{(n-2\lambda)_{k-2,\lambda}}{(k-2)!}n(n-\lambda)+(1-\lambda)n\sum_{k=1}^{\infty}\frac{(n-\lambda)_{k-1,\lambda}}{(k-1)!}\\
\nonumber &=n(n-\lambda)\sum_{k=2}^{\infty}\binom{n-2\lambda}{k-2}_{\lambda}+n(1-\lambda)\sum_{k=1}^{\infty}\binom{n-\lambda}{k-1}_{\lambda}\\
\nonumber &=n(n-\lambda)\sum_{k=0}^{\infty}\binom{n-2\lambda}{k}_{\lambda}+n(1-\lambda)\sum_{k=0}^{\infty}\binom{n-\lambda}{k}_{\lambda}\\
\nonumber &=n(n-\lambda)(1+\lambda)^{\frac{n}{\lambda}-2}+n(1-\lambda)(1+\lambda)^{\frac{n}{\lambda}-1} \\
\nonumber &=n(n+1-\lambda-\lambda^2)(1+\lambda)^{\frac{n}{\lambda}-2}.
\end{align}
By \eqref{34}, we get
\begin{equation}
	\label{35} H_{\lambda}(n,2)=n(n+1-2\lambda)(1+\lambda)^{\frac{n}{\lambda}-2}.
\end{equation}
\indent Let us take $m=3$. Then we have
\begin{align*}
	H_{\lambda}(n,3)&=\sum_{k=0}^{\infty}\binom{n}{k}_{\lambda}(k)_{3,\lambda}=\sum_{k=0}^{\infty}\binom{n}{k}_{\lambda}k(k-\lambda)(k-2\lambda) \\
	&=\sum_{k=1}^{\infty}\frac{(n)_{k,\lambda}}{(k-1)!}(k-\lambda)(k-2\lambda)=\sum_{k=1}^{\infty}\frac{(n)_{k,\lambda}}{(k-1)!}(k-1+1-\lambda)(k-2\lambda)\\
	&=\sum_{k=2}^{\infty} \frac{(n)_{k,\lambda}}{(k-2)!}(k-2\lambda)+(1-\lambda)\sum_{k=1}^{\infty}\frac{(n)_{k,\lambda}}{(k-1)!}(k-2\lambda)\\
	& =\sum_{k=3}^{\infty} \frac{(n)_{k,\lambda}}{(k-3)!}+3(1-\lambda)\sum_{k=2}^{\infty}\frac{(n)_{k,\lambda}}{(k-2)!}+(1-\lambda)(1-2\lambda)\sum_{k=1}^{\infty}\frac{(n)_{k,\lambda}}{(k-1)!} \\
	&=(n)_{3,\lambda}\sum_{k=0}^{\infty}\binom{n-3\lambda}{k}_{\lambda}+3(1-\lambda)(n)_{2,\lambda}\sum_{k=0}^{\infty}\binom{n-2\lambda}{k}_{\lambda}\\
	&\quad +(1-\lambda)(1-2\lambda)(n)_{1,\lambda}\sum_{k=0}^{\infty}\binom{n-\lambda}{k}_{\lambda}\\
	&=(n)_{3,\lambda}(1+\lambda)^{\frac{n}{\lambda}-3}+3(1-\lambda)(n)_{2,\lambda}(1+\lambda)^{\frac{n}{\lambda}-2}+(1-\lambda)(1-2\lambda)n(1+\lambda)^{\frac{n}{\lambda}-1}.
\end{align*}
Note that
\begin{displaymath}
	\lim_{\lambda\rightarrow 1}H_{\lambda}(n,3)=\sum_{k=0}^{n}\binom{n}{k}(k)_{3}=(n)_{3}2^{n-3}
\end{displaymath}
and
\begin{displaymath}
	\lim_{\lambda\rightarrow 0}H _{\lambda}(n,3)=\sum_{k=0}^{\infty}\frac{n^{k}}{k!}k^{3}=n(n^{2}+3n+1)e^{n}.
\end{displaymath}
\indent For $s\in\mathbb{C}$ with $R(s)>0$, the gamma function is defined by
\begin{displaymath}
	\Gamma(s)=\int_{0}^{\infty}e^{-t}t^{s-1}dt.
\end{displaymath}
Let $n$ be a nonnegative integer. Then
\begin{displaymath}
	\frac{\langle b\rangle_{n,\lambda}}{\langle c\rangle_{n,\lambda}}=\frac{\Gamma\big(\frac{b}{\lambda}+n\big)\Gamma\big(\frac{c}{\lambda}\big)}{\Gamma\big(\frac{c}{\lambda}+n\big)\Gamma\big(\frac{b}{\lambda}\big)},
\end{displaymath}
where $R\big(\frac{c}{\lambda}\big)>0$ and $R\big(\frac{b}{\lambda}\big)>0$. \\
For $R\big(\frac{c}{\lambda}\big)>R\big(\frac{b}{\lambda}\big)>0$, we have
\begin{equation}
	\label{36} \frac{\Gamma\big(\frac{b}{\lambda}+n\big)\Gamma\big(\frac{c}{\lambda}-\frac{b}{\lambda}\big)}{\Gamma\big(\frac{c}{\lambda}+n\big)}=\int_{0}^{1}t^{\frac{b}{\lambda}+n-1}(1-t)^{\frac{c}{\lambda}-\frac{b}{\lambda}-1}dt.
\end{equation}
From \eqref{36}, we note that
\begin{align}
	\label{37}\pFq{}{\lambda}{a,b}{c}{z}&=\sum_{n=0}^{\infty}\frac{\langle a\rangle_{n,\lambda}\langle b\rangle_{n,\lambda}}{\langle c\rangle_{n,\lambda}}\frac{z^{n}}{n!}\\ \nonumber &=\frac{\Gamma\big(\frac{c}{\lambda}\big)}{\Gamma\big(\frac{b}{\lambda}\big)\Gamma\big(\frac{c}{\lambda}-\frac{b}{\lambda}\big)}\int_{0}^{1}t^{\frac{b}{\lambda}-1}(1-t)^{\frac{c}{\lambda}-\frac{b}{\lambda}-1}(1-\lambda tz)^{-\frac{a}{\lambda}}dt.
\end{align}
In particular, for $z=\frac{1}{\lambda}$, $(\lambda\ne 0)$, from \eqref{11} we get
\begin{equation}
	\label{38} \pFq{}{\lambda}{a,b}{c}{\frac{1}{\lambda}}=\frac{\Gamma\big(\frac{c}{\lambda}\big)\Gamma\big(\frac{c}{\lambda}-\frac{b}{\lambda}-\frac{a}{\lambda}\big)}{\Gamma\big(\frac{c}{\lambda}-\frac{b}{\lambda}\big)\Gamma\big(\frac{c}{\lambda}-\frac{a}{\lambda}\big)},\quad\mathrm{where}\quad R\bigg(\frac{c}{\lambda}-\frac{b}{\lambda}-\frac{a}{\lambda}\bigg)>0.
\end{equation}
For $n\in\mathbb{N}$, by \eqref{38}, we get
\begin{equation}
	\label{39} \pFq{}{\lambda}{-n,-n}{\lambda}{\frac{1}{\lambda}}=\frac{\Gamma(1)\Gamma\big(1+\frac{2n}{\lambda}\big)}{\Gamma\big(1+\frac{n}{\lambda}\big)\Gamma\big(1+\frac{n}{\lambda}\big)}=\frac{2\lambda}{n}\frac{\Gamma\big(\frac{2n}{\lambda}\big)}{\big(\Gamma\big(\frac{n}{\lambda}\big)\big)^{2}},
\end{equation}
where $\lambda$ is a positive real number. \\
On the other hand,
\begin{equation}
	\label{40} \pFq{}{\lambda}{-n,-n}{\lambda}{\frac{1}{\lambda}}=\sum_{k=0}^{\infty}\frac{\langle -n\rangle_{k,\lambda}\langle -n\rangle_{k,\lambda}}{\langle\lambda\rangle_{k,\lambda}}\frac{\big(\frac{1}{\lambda}\big)^{k}}{k!}=\sum_{k=0}^{\infty}\lambda^{-2k}\binom{n}{k}_{\lambda}^{2}.
	\end{equation}
	Therefore, by \eqref{39} and \eqref{40}, we obtain the following theorem.
\begin{theorem}
	Let $\lambda$ be a positive real number. For $n\in\mathbb{N}$, we have
	\begin{displaymath}
		\sum_{k=0}^{\infty}\lambda^{-2k}\binom{n}{k}_{\lambda}^{2}=\frac{2\lambda}{n}\frac{\Gamma\big(\frac{2n}{\lambda}\big)}{\big(\Gamma\big(\frac{n}{\lambda}\big)\big)^{2}}.
	\end{displaymath}
\end{theorem}
\noindent Note that
\begin{displaymath}
	\lim_{\lambda\rightarrow 1}\sum_{k=0}^{\infty}\lambda^{-2k}\binom{n}{k}^{2}_{\lambda}=\sum_{k=0}^{n}\binom{n}{k}^{2}=\frac{(2n)!}{(n!)^{2}}=\binom{2n}{n}.
\end{displaymath}
\indent Now, we define the {\it{degenerate generalized hypergeometric function}} as
\begin{equation}
\label{41}\pFq{(p,q)}{\lambda}{a_{1},a_{2},\dots,a_{p}}{b_{1},\dots,b_{q}}{z}=\sum_{k=0}^{\infty}\frac{\langle a_{1}\rangle_{k,\lambda}\cdots\langle a_{p}\rangle_{k,\lambda}}{\langle b_{1}\rangle_{k,\lambda}\cdots\langle b_{q}\rangle_{k,\lambda}}\frac{z^{k}}{k!},\quad\mathrm{where}\quad |z|<1.
\end{equation}
Let $n$ be a positive integer. Then we define the {\it{degenerate hypergeometric numbers of order $p$}} by
\begin{equation}
\label{42} \pFq{(p,p-1)}{\lambda}{-n-n,\dots,-n}{\lambda,\dots,\lambda}{(-1)^{p}\lambda^{p-1}e_{\lambda}(t)}=\sum_{m=0}^{\infty}H_{\lambda}^{(p)}(n,m)\frac{t^{m}}{m!}.
\end{equation}
From \eqref{23} and \eqref{24}, we note that $H_{\lambda}(n,m)=H_{\lambda}^{(1)}(n,m)$, and $Q_{\lambda}(m,2)=H_{\lambda}^{(2)}(n,m)$. \\
\indent In \eqref{41}, we note that
\begin{align}
\label{43}\pFq{(p,p-1)}{\lambda}{-n,-n,\dots,-n}{\lambda,\dots,\lambda}{(-1)^{p}\lambda^{p-1}e_{\lambda}(t)}&=\sum_{k=0}^{\infty}\frac{\langle-n\rangle_{k,\lambda}\langle-n\rangle_{k,\lambda}\cdots\langle-n\rangle_{k,\lambda}}{\langle\lambda\rangle_{k,\lambda}\langle\lambda\rangle_{k,\lambda}\cdots\langle\lambda\rangle_{k,\lambda}}\frac{(-1)^{pk}\lambda^{(p-1)k}}{k!}e_{\lambda}^{k}(t)\\
&\nonumber=\sum_{k=0}^{\infty}\frac{\big((n)_{k,\lambda}\big)^{p}}{(k!)^{p-1}}\frac{e_{\lambda}^{k}(t)}{k!}=\sum_{k=0}^{\infty}\binom{n}{k}_{\lambda}^{p}e_{\lambda}^{k}(t) \\
&\nonumber=\sum_{m=0}^{\infty}\bigg(\sum_{k=0}^{\infty}\binom{n}{k}_{\lambda}^{p}(k)_{m,\lambda}\bigg)\frac{t^{m}}{m!}.
\end{align}
Therefore, by \eqref{42} and \eqref{43}, we obtain the following theorem.
\begin{theorem}
	For $n,p\in\mathbb{N}$, and $m\in\mathbb{N}\cup\{0\}$, we have
	\begin{displaymath}
	H_{\lambda}^{(p)}(n,m)=\sum_{k=0}^{\infty}\binom{n}{k}_{\lambda}^{p}(k)_{m,\lambda}.
	\end{displaymath}
\end{theorem}
\noindent Note that
\begin{displaymath}
	\lim_{\lambda\rightarrow 1}H_{\lambda}^{(p)}(n,m)=\sum_{k=0}^{n}\binom{n}{k}^{p}(k)_{m},\quad\mathrm{and}\quad\lim_{\lambda\rightarrow 0}H_{\lambda}^{(p)}(n,m)=\sum_{k=0}^{\infty}\frac{n^{kp}}{(k!)^{p}}k^{m}.
\end{displaymath}
\indent From \eqref{43}, we note that
\begin{equation}
	\label{44}\pFq{(p,p-1)}{\lambda}{-n,-n,\dots,-n}{\lambda,\dots,\lambda}{(-1)^{p}\lambda^{p-1}e_{\lambda}(t)}=\sum_{k=0}^{\infty}\binom{n}{k}_{\lambda}^{p}e_{\lambda}^{k}(t)
\end{equation}
\begin{align*}
	&=\sum_{k=0}^{\infty}\binom{n}{k}_{\lambda}^{p}\sum_{l=0}^{k}\binom{k}{l}\big(e_{\lambda}(t)-1\big)^{l}=\sum_{k=0}^{\infty}\binom{n}{k}_{\lambda}^{p}\sum_{l=0}^{k}(k)_{l}\sum_{m=l}^{\infty}S_{2,\lambda}(m,l)\frac{t^{m}}{m!} \\
	&=\sum_{m=0}^{\infty}\bigg(\sum_{k=0}^{\infty}\sum_{l=0}^{k}\binom{n}{k}_{\lambda}^{p}(k)_{l}S_{2,\lambda}(m,l)\bigg)\frac{t^{m}}{m!}.
\end{align*}
Therefore, by \eqref{42} and \eqref{44}, we obtain the following theorem.
\begin{theorem}
	For $n,p\in\mathbb{N}$ and $m\in\mathbb{N}\cup\{0\}$, we have
	\begin{displaymath}
		H_{\lambda}^{(p)}(n,m)=\sum_{k=0}^{\infty}\sum_{l=0}^{k}\binom{n}{k}_{\lambda}^{p}(k)_{l}S_{2,\lambda}(m,l).
	\end{displaymath}
\end{theorem}
\indent By \eqref{42} and \eqref{43}, we get
\begin{displaymath}
	H_{\lambda}^{(p)}(n,0)=\sum_{k=0}^{\infty}\binom{n}{k}^{p}_{\lambda}=\pFq{(p,p-1)}{\lambda}{-n,-n,\dots,-n}{\lambda,\dots,\lambda}{(-1)^{p}\lambda^{p-1}}
\end{displaymath}
Note that
\begin{displaymath}
	\sum_{k=0}^{n}\binom{n}{k}^{p}=\lim_{\lambda\rightarrow 1}H_{\lambda}^{(p)}(n,0)=\pFq{(p,p-1)}{}{-n,-n,\dots,-n}{1,\dots,1}{(-1)^{p}}.
\end{displaymath}
From \eqref{41}, we have
\begin{align}
	\label{45} &\pFq{(p,p-1)}{\lambda}{-n,-n,\dots,-n}{\lambda,\dots,\lambda}{(-\lambda)^{p-1}e_{\lambda}(t)}\\
	\nonumber&=\sum_{k=0}^{\infty}\binom{n}{k}_{\lambda}^{p}(-1)^{k}e_{\lambda}^{k}(t)=\sum_{m=0}^{\infty}\bigg(\sum_{k=0}^{\infty}\binom{n}{k}_{\lambda}^{p}(-1)^{k}(k)_{m,\lambda}\bigg)\frac{t^{m}}{m!}
\end{align}
Thus, by \eqref{45}, we get
\begin{equation}
	\label{46} \pFq{(p,p-1)}{\lambda}{-n,-n,\dots,-n}{\lambda,\dots,\lambda}{(-\lambda)^{p-1}}=\sum_{k=0}^{\infty}\binom{n}{k}_{\lambda}^{p}(-1)^{k}.
\end{equation}
Note that
\begin{displaymath}
	\lim_{\lambda\rightarrow 1}\pFq{(p,p-1)}{\lambda}{-n,-n,\dots,-n}{\lambda,\dots,\lambda}{(-\lambda)^{p-1}}=\sum_{k=0}^{n}\binom{n}{k}^{p}(-1)^{k}=\pFq{(p,p-1)}{}{-n,-n,\dots,-n}{1,\dots,1}{(-1)^{p-1}}.
\end{displaymath}

\section{Further Remarks.}
\indent Let $n$ be a positive integer. From \eqref{10}, we have
\begin{equation}
	\label{47} \pFq{(2,1)}{}{-n,-n}{1}{e_{\lambda}(t)}=\sum_{k=0}^{\infty}\frac{\langle-n\rangle_{k}\langle -n\rangle_{k}}{\langle 1\rangle_{k}}\frac{e_{\lambda}^{k}(t)}{k!}=\sum_{k=0}^{n}\binom{n}{k}^{2}e_{\lambda}^{k}(t)=\sum_{m=0}^{\infty}\bigg(\sum_{k=0}^{n}\binom{n}{k}^{2}(k)_{m,\lambda}\bigg)\frac{t^{m}}{m!}.
\end{equation}
Now, we define the {\it{$\lambda$-hypergeometric numbers of order $p$}} by
\begin{equation}
\label{48}	\pFq{(p,p-1)}{}{-n,-n,\dots,-n}{1,\dots,1}{(-1)^{p}e_{\lambda}(t)}=\sum_{m=0}^{\infty}H_{m,\lambda}^{(p)}(n)\frac{t^{m}}{m!}.
\end{equation}
By \eqref{47} and \eqref{48}, we get
\begin{equation}
\label{49} H_{m,\lambda}^{(2)}(n)=\sum_{k=0}^{n}\binom{n}{k}^2(k)_{m,\lambda},
\end{equation}
where $n\in\mathbb{N}$ and $m\in\mathbb{N}\cup\{0\}$. \\
~~\\
\indent The {\it{alternating $\lambda$-hypergeometric numbers of order $p$}} are defined by
\begin{equation}
	\label{50} \pFq{(p,p-1)}{}{-n,-n,\dots,-n}{1,\dots,1}{(-1)^{p-1}e_{\lambda}(t)}=\sum_{m=0}^{\infty}T_{m,\lambda}^{(p)}(n)\frac{t^{m}}{m!}.
\end{equation}
By \eqref{10}, we get
\begin{equation}
	\label{51} \pFq{(2,1)}{}{-n,-n}{1}{-e_{\lambda}(t)}=\sum_{m=0}^{\infty}\bigg(\sum_{k=0}^{\infty}\binom{n}{k}^{2}(-1)^{k}(k)_{m,\lambda}\bigg)\frac{t^{m}}{m!}.
\end{equation}
From \eqref{50} and \eqref{51}, we have
\begin{equation}
	\label{52} T_{m,\lambda}^{(2)}(n)=\sum_{k=0}^{n}\binom{n}{k}^{2}(-1)^{k}(k)_{m,\lambda},
\end{equation}
where $m\in\mathbb{N}\cup\{0\}$ and $n\in\mathbb{N}$. \\
~~\\
In general, we have
\begin{displaymath}
	T_{m,\lambda}^{(p)}(n)=\sum_{k=0}^{n}\binom{n}{k}^{p}(-1)^{k}(k)_{m,\lambda},\quad\mathrm{and}\quad H_{m,\lambda}^{(p)}(n)=\sum_{k=0}^{n}\binom{n}{k}^{p}(k)_{m,\lambda},
\end{displaymath}
where $n,p\in\mathbb{N}$ and $m\in\mathbb{N}\cup\{0\}$. \\
\indent We observe that
\begin{equation}
\label{53}	\pFq{(1,0)}{}{-n}{\times}{-e_{\lambda}(t)}=\sum_{k=0}^{n}\binom{n}{k}e_{\lambda}^{k}(t),
\end{equation}
and
\begin{equation}
	\label{54} \pFq{(1,0)}{}{-n}{\times}{e_{\lambda}(t)}=\sum_{k=0}^{n}\binom{n}{k}(-1)^{k}e_{\lambda}^{k}(t)
\end{equation}
Thus, we note that
\begin{displaymath}
	H_{m,\lambda}^{(1)}(n)=\sum_{k=0}^{n}\binom{n}{k}(k)_{m,\lambda},\quad T_{m,\lambda}^{(1)}=\sum_{k=0}^{n}\binom{n}{k}(-1)^{k}(k)_{m,\lambda},
\end{displaymath}
where $m\in\mathbb{N}\cup\{0\}$ and $n\in\mathbb{N}$. \\
\indent For example,
\begin{align*}
	H_{0,\lambda}^{(1)}(n)&=2^{n},\quad
	 H_{1,\lambda}^{(1)}(n)=n2^{n-1}\quad
	 H_{2,\lambda}^{(1)}(n)=n(n+1-2\lambda)2^{n-2}\\
	 H_{3,\lambda}^{(1)}(n)&=(n)_{3}2^{n-3}+3(n)_{2}(1-\lambda)2^{n-2}+n(1-\lambda)_{2,\lambda}2^{n-1},\dots,\\
	 T_{0,\lambda}^{(1)}(n)&=0,\quad T_{1,\lambda}^{(1)}(n)=-\delta_{1,n}, \\
	 T_{2,\lambda}^{(1)}(n)&=n(n-1)\delta_{2,n}+(\lambda-1)\delta_{n,1},\ \cdots,
\end{align*}
where $\delta_{n,k}$ is Kronecker's symbol. \\
~~\\
\indent From \eqref{50}, we note that
\begin{align}
	\label{55}\sum_{m=0}^{\infty}T_{m,\lambda}^{(1)}(n)\frac{t^{m}}{m!}&=\pFq{(1,0)}{}{-n}{\times}{e_{\lambda}(t)}=(1-e_{\lambda}(t))^{n}=(-1)^{n}\frac{n!}{n!}(e_{\lambda}(t)-1)^{n}\\
	\nonumber &=(-1)^{n}n!\sum_{m=n}^{\infty}S_{2,\lambda}(m,n)\frac{t^{m}}{m!}.
\end{align}
On the other hand,
\begin{equation}
	\label{56} (1-e_{\lambda}(t))^{n}=\sum_{j=0}^{n}\binom{n}{j}(-1)^{n-j}e_{\lambda}^{j}(t)=\sum_{m=0}^{\infty}\bigg(\sum_{j=0}^{n}\binom{n}{j}(-1)^{n-j}(j)_{m,\lambda}\bigg)\frac{t^{m}}{m!}.
\end{equation}
Thus, by \eqref{55} and \eqref{56}, we get
\begin{equation}
	\label{57} T_{m,\lambda}^{(1)}(n)=\sum_{j=0}^{n}\binom{n}{j}(-1)^{n-j}(j)_{m,\lambda}=(-1)^{n}n!S_{2,\lambda}(m,n),\quad (m\ge n).
\end{equation}
If $m<n$, then $T_{m,\lambda}^{(1)}(n)=0$.
\begin{theorem}
	For $n\in\mathbb{N}$ and $m\in\mathbb{N}\cup\{0\}$, we have
	\begin{displaymath}
		T_{m,\lambda}^{(1)}(n)=\sum_{j=0}^{n}\binom{n}{j}(-1)^{n-j}(j)_{m,\lambda}=(-1)^{n}n!S_{2,\lambda}(m,n),\quad (m\ge n).
	\end{displaymath}
	In particular, if $m<n$, then
	\begin{displaymath}
		T_{m,\lambda}^{(1)}=0.
	\end{displaymath}
\end{theorem}
For $k\ge 0$, we have
\begin{displaymath}
	\lim_{\lambda\rightarrow 0}(-1)^{n}T_{n+k,\lambda}^{(1)}(n)=\sum_{j=0}^{n}\binom{n}{j}(-1)^{j}j^{n+k}=n!S_{2}(n+k,n).
\end{displaymath}
\begin{corollary}
	For $k\ge 0$ and $n\in\mathbb{N}$, we have
	\begin{displaymath}
		\frac{1}{n!}\sum_{j=0}^{n}\binom{n}{j}(-1)^{j}j^{n+k}=S_{2}(n+k,n).
	\end{displaymath}
\end{corollary}
It is easy to show that
\begin{align*}
	(e^{t}-1)^{n} &=\sum_{j=0}^{n}\binom{n}{j}(-1)^{n-j}e^{jt}\\
	&=\sum_{j=0}^{n}\sum_{l=0}^{j}\binom{n}{j}(-1)^{n-j}(j)_{l}\frac{1}{l!}(e^{t}-1)^{l} \\
	&=\sum_{j=0}^{n}\sum_{l=0}^{j}\binom{n}{j}(j)_{l}(-1)^{n-j}\sum_{m=l}^{\infty}S_{2}(m,l)\frac{t^{m}}{m!}.
\end{align*}
Thus, we have
\begin{equation}
	\label{58} \sum_{m=n}^{\infty}S_{2}(m,n)\frac{t^{m}}{m!}=\frac{1}{n!}(e^{t}-1)^{n} =\frac{1}{n!}\sum_{j=0}^{n}\sum_{l=0}^{j}\binom{n}{j}(j)_{l}(-1)^{n-j}\sum_{m=l}^{\infty}S_{2}(m,l)\frac{t^{m}}{m!}.
\end{equation}
For $k\ge 0,\ n\in\mathbb{N}$, by \eqref{58}, we get
\begin{equation}
\label{59} S_{2}(n+k,n)=\frac{1}{n!}\sum_{j=0}^{n}\sum_{l=0}^{n+k}\binom{n}{j}(j)_{l}(-1)^{n-j}S_{2}(n+k,l).
\end{equation}
From Corollary 3.2 and \eqref{59}, we have
\begin{displaymath}
	\sum_{j=0}^{n}\binom{n}{j}(-1)^{j}j^{n+k}=\sum_{j=0}^{n}\sum_{l=0}^{n+k}\binom{n}{j}(j)_{l}(-1)^{n-j}S_{2}(n+k,l).
\end{displaymath}
\indent For $\lambda,\lambda_{1}\in\mathbb{R}$, let us define {\it{Apostol type alternating $\lambda$-hypergeometric numbers of order $p$}} by
\begin{equation}
	\label{60}\pFq{(p,p-1)}{}{-n,-n,\dots,-n}{1,\dots,1}{(-1)^{p-1}\lambda_{1}e_{\lambda}(t)}=\sum_{m=0}^{\infty}T_{m,\lambda}^{(p)}(n\ |\ \lambda_{1})\frac{t^{m}}{m!}.
\end{equation}
By \eqref{10}, we get
\begin{align}
	\label{61} \pFq{(1,0)}{}{-n}{\times}{\lambda_{1}e_{\lambda}(t)}&=\sum_{k=0}^{n}\binom{n}{k}\big(-\lambda_{1}e_{\lambda}(t)\big)^{k}=\sum_{j=0}^{n}\binom{n}{j}(-1)^{j}\lambda_{1}^{j}e_{\lambda}^{j}(t)\\
	\nonumber &=\sum_{m=0}^{\infty}\bigg(\sum_{j=0}^{n}\binom{n}{j}(-1)^{j}\lambda_{1}^{j}(j)_{m,\lambda}\bigg)\frac{t^{m}}{m!}.
\end{align}
On the other hand,
\begin{align}
	\label{62} \pFq{(1,0)}{}{-n}{\times}{\lambda_{1}e_{\lambda}(t)}&=\big(1-\lambda_{1}e_{\lambda}(t)\big)^{n}=\sum_{j=0}^{n}\binom{n}{j}(-1)^{j}\lambda_{1}^{j}e_{\lambda}^{j}(t)\\
	\nonumber &=\sum_{j=0}^{n}\binom{n}{j}(-1)^{j}\lambda_{1}^{j}\sum_{l=0}^{j}\binom{j}{l}\big(e_{\lambda}(t)-1\big)^{l}\\
	\nonumber &=\sum_{j=0}^{n}\sum_{l=0}^{j}\binom{n}{j}(-1)^{j}\lambda_{1}^{j}(j)_{l}\frac{\big(e_{\lambda}(t)-1\big)^{l}}{l!}\\
	\nonumber &=\sum_{j=0}^{n}\sum_{l=0}^{j}\binom{n}{j}(-1)^{j}\lambda_{1}^{j}(j)_{l}\sum_{m=l}^{\infty}S_{2,\lambda}(m,l)\frac{t^{m}}{m!}.
\end{align}
For $k\ge 0$ and $n\in\mathbb{N}$, by \eqref{60}, \eqref{61} and \eqref{62}, we obtain the following theorem.
\begin{theorem}
	For $\lambda,\lambda_{1}\in\mathbb{R},\ n\in\mathbb{N}$ and $k\in\mathbb{N}\cup\{0\}$, we have
	\begin{displaymath}
		T_{n+k,\lambda}^{(1)}(n~|~\lambda_{1})=\sum_{j=0}^{n}\binom{n}{j}(-1)^{j}\lambda_{1}^{j}(j)_{n+k,\lambda}=\sum_{j=0}^{n}\sum_{l=0}^{j}\binom{n}{j}(-1)^{j}\lambda_{1}^{j}(j)_{l}S_{2,\lambda}(n+k,l).
	\end{displaymath}
\end{theorem}
\noindent Note that
\begin{displaymath}
\lim_{\lambda\rightarrow 0}T_{n+k,\lambda}^{(1)}(n~|~\lambda_{1})=\sum_{j=0}^{n}\binom{n}{j}(-1)^{j}\lambda_{1}^{j}j^{n+k}=\sum_{j=0}^{n}\sum_{l=0}^{j}\binom{n}{j}(-1)^{j}\lambda_{1}^{j}(j)_{l}S_{2}(n+k,l).
	\end{displaymath}
\indent For $\lambda_{1}\in\mathbb{R}$, let us define {\it{Apstol-Stirling numbers of the second kind}} as
\begin{equation}
	\label{64}\frac{1}{n!}\big(\lambda_{1}e^{t}-1\big)^{n}=\sum_{m=0}^{\infty}S(m,n~|~\lambda_{1})\frac{t^{m}}{m!},\quad (n\ge 0).
\end{equation}
Now, we observe that
\begin{equation}
	\label{65}\frac{1}{n!}(\lambda_{1}e^{t}-1)^{n}=\frac{1}{n!}\sum_{j=0}^{n}\binom{n}{j}\lambda_{1}^{j}(-1)^{n-j}e^{jt}=\frac{1}{n!}\sum_{j=0}^{n}\binom{n}{j}\lambda_{1}^{j}(-1)^{n-j}\sum_{m=0}^{\infty}j^{m}\frac{t^{m}}{m!}.
\end{equation}
For $k\in\mathbb{N}\cup\{0\}$ and $n\in\mathbb{N}$, by \eqref{64} and \eqref{65}, we get
\begin{equation}
	\label{66} S(n+k,\ n~|~\lambda_{1})=\frac{(-1)^{n}}{n!}\sum_{j=0}^{n}\binom{n}{j}\lambda_{1}^{j}(-1)^{j}j^{n+k}.
\end{equation}
From Theorem 3.3 and \eqref{66}, we have
\begin{equation}
	\label{67}\lim_{\lambda\rightarrow 0}T_{n+k,\lambda}^{(1)}(n~|~\lambda_{1})=(-1)^{n}n!S(n+k,n~|~\lambda_{1})=\sum_{j=0}^{n}\binom{n}{j}\lambda_{1}^{j}(-1)^{j}j^{n+k}.
\end{equation}
Therefore, by \eqref{67}, we obtain the following corollary.
\begin{corollary}
	For $k\in\mathbb{N}\cup\{0\}$ and $n\in\mathbb{N}$, we have
	\begin{displaymath}
		\sum_{j=0}^{n}\binom{n}{j}\lambda_{1}^{j}(-1)^{j}j^{n+k}=(-1)^{n}n!S(n+k,n~|~\lambda)=\lim_{\lambda\rightarrow 0}T_{n+k,\lambda}^{(1)}(n~|~\lambda_{1}).
	\end{displaymath}
\end{corollary}

\emph{Remarks.}
(a) The Corollary 3.4 naturally interprets the sum in Corollary 3.4 in terms of Apostol-Stirling numbers of the second kind defined by \eqref{64}. In answering an open question asked by Srivastava, Simsek gave three different expressions for this sum as special cases of the results in [17,Theorems 17-19]. \\
(b) For $\lambda,\lambda_{1}\in\mathbb{R}$, let us define {\it{Apostol type $\lambda$-hypergeometric numbers of order $p$}} by
\begin{equation}
\label{68} \pFq{(p,p-1)}{}{-n,-n,\dots,-n}{1,\dots,1}{(-1)^{p}\lambda_{1}e_{\lambda}(t)}=\sum_{m=0}^{\infty}H_{m,\lambda}^{(p)}(n~|~\lambda_{1})\frac{t^{m}}{m!}.
\end{equation}
By \eqref{10}, we get
\begin{equation}
	\label{69} \pFq{(p,p-1)}{}{-n,-n,\dots,-n}{1,\dots,1}{(-1)^{p}\lambda_{1}e_{\lambda}(t)}=\sum_{k=0}^{\infty}\binom{n}{k}^{p}\lambda_{1}^{k}e_{\lambda}^{k}(t)=\sum_{m=0}^{\infty}\bigg(\sum_{k=0}^{n}\binom{n}{k}^{p}\lambda_{1}^{k}(k)_{m,\lambda}\bigg)\frac{t^{m}}{m!}.
\end{equation}
From \eqref{68} and \eqref{69}, we have
\begin{displaymath}
	H_{m,\lambda}^{(p)}(n~|~\lambda_{1})=\sum_{k=0}^{n}\binom{n}{k}^{p}\lambda_{1}^{p}(k)_{m,\lambda},\quad\mathrm{where}\ n,p\in\mathbb{N}\ \mathrm{and}\ m\in\mathbb{N}\cup\{0\}.
\end{displaymath}
Note that
\begin{displaymath}
	\lim_{\lambda\rightarrow 0}H_{m,\lambda}^{(p)}(n~|~\lambda_1)=\sum_{k=0}^{n}\binom{n}{k}^{p}\lambda_{1}^{p}k^{m}.
\end{displaymath}


\begin{thebibliography}{9}
\bibitem{1}
L. Carlitz, \emph{Degenerate Stirling, Bernoulli and Eulerian numbers,} Utilitas Math. \textbf{15} (1979), 51-88.
\bibitem{2}
L. Carlitz, \emph{A degenerate Staudt-Clausen theorem,} Arch. Math. (Basel) \textbf{7} (1956), 28-33.
\bibitem{3}
D. V. Dolgy, T. Kim, \emph{Some explicit formulas of degenerate Stirling numbers associated with the degenerate special numbers and polynomials,} Proc. Jangjeon Math. Soc. \textbf{21} (2018), no. 2, 309-317.
\bibitem{4}
R. Golombek, \emph{Aufgabe 1088,} Elem. Math. \textbf{49} (1994), 126-127.
\bibitem{5}
 R. Golombek, D. Marburg, \emph{Aufgabe 1088, Summen mit Quadraten von Binomialkoeffizienten,} Elem. Math., \textbf{50} (1995), 125-131.
\bibitem{6}
H. Haroon, W. A. Khan, \emph{Degenerate Bernoulli numbers and polynomials associated with degenerate Hermite polynomials,} Commun. Korean Math. Soc. \textbf{33} (2018). no. 2, 651-669.
\bibitem{7}
Y. He, S. Araci, \emph{Sums of products of Apostol-Bernoulli and Apostol-Euler polynomials,} Adv. Difference Equ. (2014), 2014:155, 13 pp.
\bibitem{8}
D. S. Kim, H. Y. Kim, D. Kim, T. Kim, \emph{Identities of Symmetry for Type 2 Bernoulli and Euler Polynomials,} Symmetry 2019, \textbf{11}, 613.
\bibitem{9}
T. Kim, G.-W. Jang, \emph{A note on degenerate gamma function and degenerate Stirling number of the second kind,} Adv. Stud. Contemp. Math. (Kyungshang) \textbf{28} (2018), no. 2, 207-214.
\bibitem{10}
T. Kim, D. S. Kim, \emph{Degenerate Laplace transform and degenerate gamma function,} Russ. J. Math. Phys. 24 (2017), no. \textbf{2}, 241-248.
\bibitem{11}
T. Kim, \emph{A note on degenerate Stirling polynomials of the second kind,} Proc. Jangjeon Math. Soc. \textbf{20} (2017), no. 3, 319-331.
\bibitem{12}
T. Kim, Y. Yao, D. S. Kim, G. W. Jang, \emph{Degenerate-Stirling numbers and $r$-Bell polynomials,} Russ. J. Math. Phys. 25 (2018). no. \textbf{1}, 44-58.
\bibitem{13}
T. Kim, D. S. Kim, \emph{A note on type 2 Changhee and Daehee polynomials,} Rev. R. Acad. Cienc. Exactas Fis. Nat. Ser. A Mat. RACSAM 113 (2019), no. \textbf{3}, 2783-2791.
\bibitem{14}
S.-S. Pyo, T. Kim, \emph{Some identities of fully degenerate Bell polynomials arising from differential equations,} Proc. Jangjeon Math. Soc. 22 (2019), no. \textbf{2}, 357-363.
\bibitem{15}
E. D. Rainville, \emph{Special Functions,} Chelsea Publ. Comp., Bronx, New York, 1971.
\bibitem{16}
S. Roman, \emph{The umbral calculus,} Pure and Applied Mathematics, 111, Academic Press, Inc. [Harcourt Brace Jovanovich. Publishers], New York, 1984.
\bibitem{17}
Y. Simsek, \emph{Generating functions for finite sums involving higher powers of binomial coefficients: analysis of hypergeometric functions including new families of polynomials and numbers,} J. Math. Anal. Appl. 477 (2019), no. \textbf{2}, 1328-1352.
\bibitem{18}
Y. Simsek, \emph{Identities and relations related to combinatorial numbers and polynomials,} Proc. Jangjeon Math. Soc. 20 (2017), no. \textbf{1}, 127-135.
\bibitem{19}
Y. Simsek, \emph{Identities on the Changhee numbers and Apostol-type Daehee polynomials,} Adv. Stud. Contemp. Math.( Kyungshang) 27 (2017), no. \textbf{2}, 199-212.
\bibitem{20}
E. T. Whittaker, G. N. Watson, \emph{Modern Analysis, 4th ed.,} Cambridge: Cambridge Univ. Press, 1927.
\bibitem{21}
W. Zhang, X. Lin, \emph{Identities involving trigonometric functions and Bernoulli numbers,} Appl. Math. Comput. 334 (2018), 288-294.
\end{thebibliography}
\end{document}